 \documentclass[journal,12pt,onecolumn,draftclsnofoot]{IEEEtran}  %

\IEEEoverridecommandlockouts                              

\usepackage{amsmath,amssymb, graphicx,epstopdf,textcomp}
\usepackage{enumerate,url,algpseudocode}
\usepackage[boxruled]{algorithm2e}
\usepackage{amsthm}
\allowdisplaybreaks
\usepackage{cite}
\usepackage[tight,footnotesize]{subfigure}
\usepackage[usenames,dvipsnames]{color}
\usepackage{dblfloatfix}
\usepackage{tablefootnote,booktabs}

\newcommand{\IEEEbq}{\begin{IEEEeqnarray}{rCl}}
	\newcommand{\IEEEeq}{\end{IEEEeqnarray}}
\newcommand{\bsq}{\begin{subequations}}
	\newcommand{\esq}{\end{subequations}}

\newcommand\oprocendsymbol{\hbox{$\square$}}
\newcommand\oprocend{\relax\ifmmode\else\unskip\hfill\fi\oprocendsymbol}

\def \N {\mathcal N}
\def\J{\mathcal J}
\def\R{\mathcal R}
\def\T{\mathcal T}
\def\L{\mathcal L}
\def\P{\mathcal P}
\def\V{\mathcal V}

\def\M{\mathcal M}
\def\E{\mathcal E}

\newtheorem{theorem}{Theorem}

\newtheorem{claim}[theorem]{Claim}

\newcommand{\bp}{{\bf p}}
\newcommand{\bq}{{\bf q}}

\newcommand{\revf}[1]{{#1}}
\newcommand{\revs}[1]{{#1}}

\newif\ifwidefig

\widefigfalse

\begin{document}
	
	\title{Optimal Water-Power Flow Problem: \\ Formulation and  Distributed Optimal Solution}
	\author{Ahmed S. Zamzam, Emiliano Dall'Anese,  Changhong Zhao, Josh A. Taylor,  and Nicholas D. Sidiropoulos
		\thanks{A. S. Zamzam is with the Department of Electrical and Computer Engineering, University of Minnesota, Minneapolis, MN 55455, USA. Emails: ahmedz@umn.edu}%
		\thanks{E. Dall'Anese  and C. Zhao are with the National Renewable Energy Laboratory, Golden, CO 80401, USA. Emails: \{emiliano.dallanese, changhong.zhao\}@nrel.gov}
		\thanks{J. A. Taylor is with the Department of Electrical and Computer Engineering, University of Toronto, Toronto, ON M5S 3G4, Canada. Email: josh.taylor@utoronto.ca }
		\thanks{N. D. Sidiropoulos is with the Department of Electrical and Computer Engineering, University of Virginia, Charlottesville, VA 22904, USA. Email: nikos@virginia.edu }
		\thanks{The work of A.S. Zamzam and N.D. Sidiropoulos was partially supported by NSF under grant CIF-1525194. The work of E. Dall'Anese and C. Zhao was supported in part by the Laboratory Directed Research and Development program at the National Renewable Energy Laboratory.}
	}
	
	\maketitle
	\begin{abstract}
		This paper formalizes an optimal water-power flow (OWPF) problem to optimize the use of controllable assets across power and water systems while accounting for the couplings between the two infrastructures. Tanks and pumps are optimally managed to satisfy water demand while improving power grid operations; {for the power network, an AC optimal power flow formulation is augmented to  accommodate the controllability of water pumps.} Unfortunately, the physics governing the operation of the two infrastructures and coupling constraints lead to a nonconvex (and, in fact, NP-hard) problem; however, after reformulating OWPF as a nonconvex, quadratically-constrained quadratic problem, a feasible point pursuit-successive convex approximation approach is used to identify feasible and optimal solutions. In addition, a distributed solver based on the alternating direction method of multipliers enables water and power operators to pursue individual objectives while respecting the couplings between the two networks. The merits of the proposed approach are demonstrated for the case of a distribution feeder coupled with a municipal water distribution network.
	\end{abstract}
	
	
	\section{Introduction}
	\label{sec:overview}
	
	Power and water networks are critical infrastructures. These systems are predominantly planned and operated independently, although their operation is intrinsically coupled at multiple spatial and temporal scales.  For example, electric  pumps for agricultural and municipal water systems affect the operation (in the form of power and energy demands) of power distribution grids~\cite{albadi2008summary, siano2014demand,marks2014opportunities}; on the other hand, the capacity of thermoelectric power plants strongly depends on the availability of cooling water~\cite{byers2014electricity, colman2013effect, wang2016assessing}. With reference to municipal and wastewater systems, electric pumps are key elements to  overcoming  geographical differences in head pressure and head losses caused by pipe friction, and they enable the supply of water demand within given water quality standards. Electricity consumption from pumps constitutes major operating costs for water utilities. For example, in the United States the overall  operation of drinking and wastewater networks  represents  $4\%$ of the total electricity consumption~\cite{denig2008reducing}. Optimizing water pump operation has therefore significant potential to save energy, reduce emissions, and enhance the reliability and efficiency of the power grid.
	
	Under dynamic electricity pricing, water utilities can schedule their pumps and adjust the consumption of variable-speed-drive pumps to minimize the cost of energy. The optimal pump scheduling problem is often formulated as a (nonconvex) mixed integer nonlinear program (MINLP), wherein nonlinearity stems from the water network hydraulic model, and binary variables indicate the pumps' on/off status. In the literature, various approaches have been developed for MINLPs. These include (piecewise) linear approximations \cite{jowitt1992optimal, giacomello2012fast, verleye2013optimising,sherali2001effective}, nonconvex nonlinear programming relaxation \cite{burgschweiger2009optimization}, continuous constraint relaxation combined with branch and bound \cite{ulanicki2007dynamic}, hydraulic simulation to implicitly enforce nonlinear constraints \cite{brion1991methodology}, Lagrange decomposition integrated with simulation-based search \cite{ghaddar2015lagrangian}, gradient method together with sensitivity analysis \cite{yu1994optimized}, and problem-specific presolving techniques \cite{gleixner2012towards}. Finally,~\cite{fooladivanda2017energy} bypasses the nonlinearity of the water network hydraulic model by leveraging a second-order \revs{cone} relaxation, and sufficient conditions for the exactness of the relaxation are provided.  Optimizing the operation of municipal water distribution networks  in response to time-varying electricity prices was also considered in~\cite{Oikonomou-2017}. A two-step approach is taken whereby the operation of the water network is optimized \emph{a priori}, and the controllable assets in the power network are then optimized in a subsequent stage based on the consumption of the water pumps. 
	
	A variety of additional objectives and constraints can arise when optimizing the operation of water pumps. 
	Examples include the cost of energy loss caused by pipe friction \cite{fooladivanda2017energy}, maximum electricity demand charges (together with unit charges) \cite{yu1994optimized}, peak power consumption, maintenance costs, reservoir level variation \cite{baran2005multi}, and land subsidence \cite{wang2009enhanced}. 
	Another line of work focuses on dynamic optimal pump management, which is usually formulated as a dynamic programming problem \cite{zessler1989optimal, nitivattananon1996optimization, lansey1994optimal}. A difficulty is to obtain a near optimal solution of the dynamic program in a reasonable amount of time. To overcome this challenge, techniques such as simplifying hydraulic dynamics \cite{lansey1994optimal}, spatial decomposition of water networks, and temporal decomposition of operation periods \cite{zessler1989optimal, nitivattananon1996optimization} are exploited in formulating and solving the dynamic program.
	
	
	The works mentioned above pertain to optimal management of water networks, and the main interaction with power systems\textemdash if any\textemdash is in the form of responsiveness to electricity prices. It is, however,  increasingly recognized that a \emph{joint} optimization of power and water infrastructures can bring significant benefits from operational and economical standpoints~\cite{MagDallAnese}. Controllable assets of water utilities can provide valuable services to power systems to enhance reliability and efficiency, as well as to cope with the volatility of distributed renewable generation; these services include frequency regulation, regulating reserves, and contingency reserves. On the other hand, the incentives for the provisioning of grid services to electric utilities could be used by water system operators for capital improvements and capacity expansion. To the best of our knowledge, however, there are no systematic approaches to jointly optimize the use of controllable assets across power and water systems while acknowledging intrinsic couplings between the two infrastructures. 
	

	This paper formulates an optimal water-power flow (OWPF) problem to minimize the (sum of the) cost functions associated with water and power operators while respecting relevant engineering and operational constraints of the two systems as well as pertinent intra-system coupling constraints. In particular, the power consumed by a pump is related to the pump's pressure gain and flow rate. The problem is tailored to coupled power distribution feeders and municipal water distribution networks (although  its applicability is not limited to this operational setting), and it addresses the controllability of distributed energy resources (DERs) and water pumps. Pump selection is not addressed because it is assumed that a pump scheduling problem is solved at a slower timescale than the  OWPF. Because of the AC power flow equations and  the water network hydraulic model, the OWPF problem is nonconvex (and NP-hard). 
	
	The proposed technical approach involves the reformulation of nonconvex constraints as equivalent nonconvex quadratic constraints, thus leading to an equivalent OWPF problem that is in the prototypical form of a nonconvex quadratically constrained quadratic program (QCQP). The resulting nonconvex QCQP is then solved by using the feasible point pursuit-successive convex approximation (FPP-SCA) method ~\cite{Mehanna-2015,Zamzam-2017}.  The FPP-SCA algorithm replaces the nonconvex constraints by inner convex surrogates around a specific point to construct a convex restriction of the original problem. Because such restriction might lead to infeasibility, the main operating principles of the algorithm to identify a feasible and optimal point involve the following two phases:
	\begin{enumerate}
		\item \emph{Feasibility phase}:  a feasible solution is obtained by solving a sequence of approximations of the original problem, where slack variables are added after restriction to quantify and ultimately zero-out the amount of constraint violations;
		\item \emph{Refinement phase}: successive convex approximation of the feasible set is used to find a Karush-Kuhn-Tucker (KKT) point of the OWPF problem.
	\end{enumerate}

	One of the practical challenges of the solution outlined above is the need for a centralized computational platform that can solve the OWPF problem. To bypass the need for a central controller, we develop a distributed solver based on the alternating direction method of multipliers~\cite{boyd2011distributed}. With the distributed solution method,  water and power operators can pursue individual operational objectives and retain controllability of their own assets (pumps for the water network and DERs for the power network) while respecting operational couplings between the two networks~\cite{MagDallAnese}. In the resulting iterative procedure,  each system operator solves a smaller optimization problem with  variables associated only with its network and controllable devices; water and power operators subsequently  exchange information regarding the shared variables to \emph{reach consensus} on the powers consumed by the water pumps.  This represents an additional unique contribution of the present paper. 
	
	Centralized and distributed methods are tested using an IEEE distribution test feeder connected to a municipal water distribution network adopted from the literature.
	
	It is worth pointing out that for the AC optimal power flow problem, a number of approaches have been proposed in the literature based on semidefinite relaxations,  second-order cone relaxation, and  linearization methods for the AC power flow equations; see, for example,~\cite{jabr2006radial, molzahn-2013, DallAnese13, low2014convex, Robbins16, Gan14,DallAnese17} and pertinent references therein. In this paper, we adopt the FPP-SCA method~\cite{Zamzam-2017} because it allows us to tackle the nonconvexity of constraints of both water and power networks in a unified manner.

	The rest of this paper is organized as follows. Section II contains the modeling of a municipal water network, and Section III introduces the model for the power distribution network. Section IV outlines the proposed OWPF problem is outlined. Section V illustrates the application of the FPP-SCA approach to the joint optimization problem, and Section VI presents the distributed algorithm. Section VII demonstrates the effectiveness of the proposed algorithm via a test case from the literature, and Section VIII summarizes conclusions and findings.
	
	\noindent {\bf Notation:} matrices (vectors) are denoted by boldface capital (small) letters unless otherwise stated; $(\cdot)^T$, $ (\cdot)^\ast $, and $(\cdot)^H$ stand for transpose, conjugate, and complex-conjugate transpose, respectively; and $|(\cdot)|$ denotes the magnitude of a number or the cardinality of a set.

	\section{Modeling the Water Network}
	\label{sec:water-model}
	
\begin{figure}[t]
		\centering
		\includegraphics[width=0.75\columnwidth]{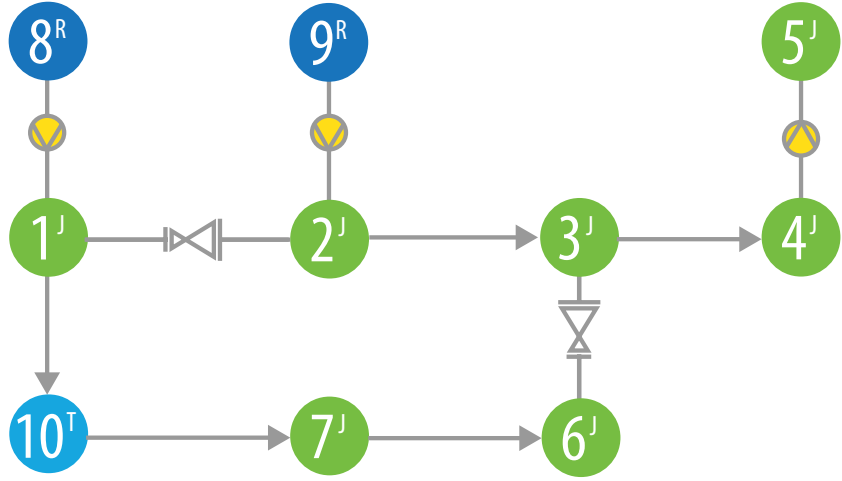}
		\caption{Example of municipal water network~\cite{Cohen-2000}; this network will be utilized in the numerical experiments.}
		\label{fig:water-network}
	\end{figure}

	As in~\cite{fooladivanda2017energy}, we consider a  municipal water network comprising a set of nodes $\mathcal{N}$ and a set of pipes $\mathcal{L}$ linking the nodes. The disjoint sets $\mathcal{J}$, $\mathcal{T}$, and $\R$, with $\N = \J \cup \T \cup \R$,  denote the sets of junctions, tanks, and reservoirs, respectively. 
	The disjoint subsets $\mathcal{P}$ and $\V$ of $\L$ have variable-speed pumps and pressure-reducing valves, respectively. We optimize the water system operation during an interval $t = 1,\dots,T$, with $\delta$ representing the time interval between two consecutive time periods. We assume that the on/off states of the pumps and valves are determined at a slower timescale and hence constant over $t=1,\dots,T$. Note that water can only flow in one direction through a pump or a valve. Moreover, we assume that the directions of water flow in pipes without pumps or valves do not change over $t=1,\dots,T$. Therefore, $(\N,\L)$ can be modeled as a \emph{directed} graph, wherein $\mathcal{P}$ and $\V$ are the sets of pipes with pumps and valves that are in an ``on'' state during $t=1,\dots,T$, respectively. 
	Let $q_{ij}^t$ denote the volumetric water flow rate through a pipe $ij \in \L$ at time $t = 1,\dots,T$. Next, we list the models for different components in water networks. \revs{Fig.~\ref{fig:water-network} illustrates an example of water network taken from~\cite{Cohen-2000}; the network features $ 7 $ junctions, $ 2 $ reservoirs, a tank, $ 3 $ pumps, and $ 2 $ valves. This network will be utilized in the numerical experiments in Section~\ref{sec:experiments}.}
	
	
	%
	%

	\emph{Junctions:} Denote the water demand at junction $j\in\mathcal{J}$ at time $t$ as $d_j^t$. Then, the following mass conservation constraint must hold:
	\begin{IEEEeqnarray}{rCl} 
		\sum_{i: ij\in\mathcal{L}}q_{ij}^t - \sum_{k: jk \in\mathcal{L}}q_{jk}^t = d_j^t, \quad \forall j\in\mathcal{J},~\forall t=1,\dots,T.\label{eq:junction-conservation}
	\end{IEEEeqnarray}
	\revf{The head pressure at junction	$j \in \mathcal{J}$ at time interval $t$, denoted by $h_j^t$, must satisfy the condition: } 
	\begin{IEEEeqnarray}{rCl}\label{eq:minimum-junction-head}
		h_j^t \geq h^{\text{min}}_j ,\quad  \forall j\in\mathcal{J},~\forall t=1,\dots,T
	\end{IEEEeqnarray}
	where $h^{\text{min}}_j$ is the minimum allowable pressure head at junction $j$ (assigned by the water system operator based on engineering considerations). 
	
	\emph{Tanks:} Let $\mu_k^t$ denote the volume of water in a tank $k \in \T$ at time $t$. The pressure head $h_k^{t,\textnormal{out}}$ at the outlet of tank $k$ at time $t$ satisfies the equation $A_k h_k^{t,\textnormal{out}} = \mu_k^t$, where the constant $A_k$ denotes the cross-sectional area of tank $k$.  
	Further, the tank outlet pressure heads satisfy the following constraints:
	\begin{IEEEeqnarray}{rCl}
		& &h_k^{t,\textnormal{out}}=h_k^{t-1,\textnormal{out}} + \frac{\delta}{A_k}  \left(\sum_{i: ik\in\mathcal{L}}q_{ik}^t - \sum_{j: kj\in\mathcal{L}}q_{kj}^t \right)  \label{eq:tank-dynamics}\\
		& &0 \leq h_k^{t,\textnormal{out}}\leq  h^{t, \text{in}}_k, \quad \forall k \in\mathcal{T},~\forall t=1,\dots, T   \label{eq:tank-h-limits}
	\end{IEEEeqnarray}
	with the initial pressure head  $h_k^{0,\textnormal{out}}$ determined based on the initial volume \revs{$\mu_k^0$} as $h_k^{0, \textnormal{out}} = \mu_k^0 / A_k$. The inlet pressure head $h^{t, \text{in}}_k$ is a variable that can be adjusted to determine the inlet flow rate. An upstream node of tank $k$ sees its inlet head $h^{t, \text{in}}_k$, whereas a downstream node sees its outlet head $h^{t, \textnormal{out}}_k$. In the sequel, for notational simplicity, we use only $h_k^t$ to denote $h^{t, \text{in}}_k$ in an equation related to an upstream node, and $h_k^{t,\textnormal{out}}$ in an equation related to a downstream node. \revs{In Fig.~\ref{fig:water-network}, the tank corresponds to node $ 10 $.}

	\emph{Reservoirs:} We treat the reservoirs $r \in \R$ as infinite sources of water where the mass conservation constraint \eqref{eq:junction-conservation} is not imposed. The pressure heads at reservoirs are set to zero: 
	\begin{IEEEeqnarray}{rCl}\label{eq:zero-head-reservoir}
		h_r^t = 0,\quad \forall r \in \R,~\forall t =1,\dots,T.
	\end{IEEEeqnarray} 
	\revs{The reservoirs in Fig.~\ref{fig:water-network} are  color-coded  in blue and are located at nodes $ 8 $ and $ 9 $.}
	
	\emph{Variable-speed pumps:} The pressure head gain due to a variable-speed pump\revs{, which is denoted by $ \hat{h}^t_{ij} $,} is modeled as\revf{~\cite{coulbeck1991ginas, ulanicki2008modeling}}:
	\begin{IEEEeqnarray}{rCl}
		&& \hat{h}^t_{ij} = A_{ij} \left(q_{ij}^t\right)^2 + B_{ij} q_{ij}^t w_{ij}^t + C_{ij} \left(w_{ij}^t\right)^2 \label{eq:pump-head-gain} \\
		&& 0 \leq  w_{ij}^t \leq w^{\text{max}}_{ij}, \quad \forall ij \in \mathcal{P}, ~\forall t=1,\dots,T \label{eq:max-speed}
	\end{IEEEeqnarray}
	where $w_{ij}^t$ and $w^{\text{max}}_{ij}$ are the ratio of the actual pump speed to the nominal speed and the maximum allowable $w_{ij}^t$, respectively. The coefficients $A_{ij} \leq 0$, $B_{ij} \geq 0$, and $C_{ij} \geq 0$ are pump parameters evaluated at the nominal speed.
	The following constraints hold for the pumps with an ``on'' state:
	\begin{IEEEeqnarray}{rCl}
		&&\hat{h}_{ij}^t =  \left(h_j^t+\overline{h}_j \right) - \left(h_i^t+\overline{h}_i\right)   \label{eq:on-pump-head-constraint}  
		\\
		&& \hat{h}_{ij}^t \geq 0    \label{eq:on-pump-head-gain-positive} 
		\\
		&& q_{ij}^t \geq  0 , \quad \forall ij \in\mathcal{P}, ~\forall t=1,\dots,T  \label{eq:on-pump-flow-rate} 
	\end{IEEEeqnarray}
	where $\overline{h}_i$ denotes the elevation of node $i$. 
	\revs{In Fig.~\ref{fig:water-network},  the pumps are located on the pipes $8 \rightarrow 1$, $9 \rightarrow 2$, and $4 \rightarrow 5$.}
	
	\emph{Pressure-reducing valves:} Let $R_{ij}^t$ denote the  head loss along a pipe $ij$ with a valve. Then, the following equations pertain to a pipe with a valve in an ``on'' state:
	\begin{IEEEeqnarray}{rCl}
		&&R_{ij}^t =  \left(h_i^t+\overline{h}_i \right) - \left(h_j^t+\overline{h}_j\right)   \label{eq:on-valve-head-loss}  \\
		&&q_{ij}^t \geq  0, \quad \forall ij \in\mathcal{V}, ~\forall t=1,\dots,T.  \label{eq:on-valve-flow-rate} 
	\end{IEEEeqnarray}
	\revs{In Fig.~\ref{fig:water-network}, two valves are present on the pipe between the junctions $ 1 $ and $ 2 $, as well as between junctions $ 3 $ and $ 6 $. }
	

	\emph{Pipes without pumps or valves:} Nodal pressure heads are related by head losses. The following constraints must be satisfied for pipes without pumps or valves:
	\begin{eqnarray}
	&&{\tilde h}_{ij}^t =  \left(h_i^t+\overline{h}_i \right) - \left(h_j^t+\overline{h}_j\right)   \label{eq:pipe-head-loss} \\
	&& q_{ij}^t \geq  0 , \quad \forall ij \in\mathcal{\L / (\P \cup \V)}, ~\forall t=1,\dots,T. \label{eq:on-pipe-flow-rate} 
	\end{eqnarray}
	The head loss, ${\tilde h}_{ij}^t$, can be approximated by the following Darcy-Weisbach equation:
	\begin{equation}
	{\tilde h}_{ij}^t = F_{ij}\left(q_{ij}^t\right)^2, \quad \forall ij \in\mathcal{\L / (\P \cup \V)}, ~\forall t=1,\dots,T  \label{eq:Darcy-Weisbach}
	\end{equation}
	where the parameter $F_{ij}$ depends on the pipe characteristics. 
	
	\emph{Power consumption of pumps:}  Let $\rho$ denote the  water density, $g$ the standard gravity coefficient, and $\eta_{ij}$ the efficiency of pump $ij$. 
	Then, the power consumption of pump $ij$ at time $t$ is given by the following expression\revf{~\cite{ulanicki2008modeling}}:
	\begin{IEEEeqnarray}{rCl}\label{eq:pump-power-consumption}
		p_{\textrm{pump},ij}^t = \frac{\rho g}{\eta_{ij}} \hat{h}_{ij}^t q_{ij}^t
	\end{IEEEeqnarray}
	for all $t=1,\dots,T. $
	
	Equation~\eqref{eq:pump-power-consumption} captures  a fundamental \emph{coupling} between the water network and the power system because it relates the electrical power consumption of a pump with the volumetric water flow rate and the head gain.

	\section{Modeling Power Distribution Networks}\label{sec:power-model}
	
	In this section, we outline the so-called bus injection model for a distribution feeder~\cite{low2014convex}. The advantage of the bus injection model is that it facilitates the formulation of the AC power flow equations in a (nonconvex) quadratic form~\cite{Zamzam-2017}. For ease of exposition and notational simplicity, the model is outlined for a balanced system; however, the proposed technical approach naturally extends to unbalanced multiphase systems as shown in~\cite{Zamzam-2017}.  
	
	\begin{figure}[t]
		\centering
		\includegraphics[width=0.6\columnwidth]{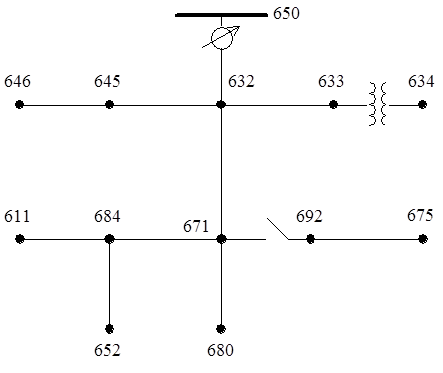}
		\caption{Example of distribution network: the IEEE 13-node feeder. This feeder will be utilized in the experiments.}
		\label{fig:power-network}
	\end{figure}

	Let $\M = \{0, 1,\dots,M\}$ denote the set of  nodes of the power distribution network. Assume that node $0$ corresponds to  the substation, and define $\M^+ := \M \setminus \{0\}$. Let $\E$ denote the set of distribution lines connecting the nodes. \revs{An example of power distribution network is illustrated in Fig.~\ref{fig:power-network}; the distribution network  is formed by $ 13 $ buses and $ 12 $ lines;  the substation (i.e., the point of connection with the transmission grid) is located at node $ 650 $. }
	For each node $m \in \M$, let $v_m^t$ denote its complex voltage, and let $s_m^t$ denote the net complex power injection at the same node at time $t = 1,\dots, T$. For each line $(m,n) \in \E$, let $z_{mn}$ and $ y_{mn} $ denote the complex series impedance and shunt admittance of a $ \pi $-equivalent circuit model. The admittance matrix $ {\bf Y} $ is obtained by setting to  $ -z_{mn}^{-1} $, $(m,n) \in \E$,  its off-diagonal elements, whereas the $m $-th diagonal element is given by $ \sum_{n: (m,n) \in \E} (z_{mn}^{-1}+\frac{1}{2} y_{mn}) $. Accordingly, using Ohm's Law and Kirchhoff's Law, the vector of current injections ${\bf i}^t := [i_0^t, i_1^t, \ldots, i_M^t]^T$ and the vector of voltages ${\bf v}^t := [v_0^t, v_1^t, \ldots, v_M^t]^T$ are related as: 
	\begin{equation}\label{nodal-currents}
	{\bf i}^t = {\bf Y} {\bf v}^t
	\end{equation}
	The net complex power injection $s_m^t $ is given by $s_m^t = v_m^t (i_m^t)^\star $. Through standard manipulations, the net active and reactive power injections at node $m$ can be expressed as:
	\begin{eqnarray}\label{power-real}
	p_m^t = {{\bf v}^t}^H {\bf Y}_m {\bf v}^t\\
	q_m^t = {{\bf v}^t}^H \overline{\bf Y}_m {\bf v}^t\label{power-imag}
	\end{eqnarray}
	where the Hermitian matrices $ {\bf Y}_m $ and $ \overline{\bf Y}_m $ are given by: 
	\begin{eqnarray}
	{\bf Y}_m = \frac{1}{2} ({\bf e}_m {\bf e}_m^T {\bf Y} + {\bf Y}^H {\bf e}_m {\bf e}_m^T)\\
	\overline{\bf Y}_m = \frac{j}{2} ({\bf e}_m {\bf e}_m^T {\bf Y} - {\bf Y}^H {\bf e}_m {\bf e}_m^T)
	\end{eqnarray}
	and  $ {\bf e}_m $ represents the $ (m+1) $-th basis of $ \mathbb{R}^{M+1}$.  Similarly, the squared voltage magnitude at node $ m $ can be expressed as:
	\begin{equation}\label{voltage-magnitude}
	|v_m^t|^2 = {{\bf v}^t}^H {\bf M}_m {\bf v}^t
	\end{equation}
	where  ${\bf M}_m := {\bf e}_m {\bf e}_m^T$.

	Define $ \mathcal{D} \subseteq {\M^+} $ to be the subset of nodes where controllable DERs are located, and let $ p_{r,m}^t $ and $ q_{r,m}^t $ denote the active and reactive power generations from the DER(s) at node $ m $ at time $ t $. The power setpoints for a DER at node $m$ are assumed to lie within a set: 
	\begin{equation}\label{setpoints}
	(p_{r,m}^t, q_{r,m}^t) \in \mathcal{C}_m^t \subset \mathbb{R}^2
	\end{equation}
	which is assumed convex. For DERs such as photovoltaic (PV) systems, energy storage systems, and electric vehicles,  the set $ \mathcal{C}_m^t $ is described by linear inequality constraints or convex quadratic constraints~\cite{Zamzam-2017,li2011optimal,low2014convex}.  Let $ p_{L,m}^t $ and $ q_{L,m}^t $ denote the (noncontrollable) active and reactive loads at node $ m \in \M^+ $ at time $ t $. Then, the net real and reactive  power injections at node $ m \in \M^+ $ are given by: 
	\begin{eqnarray}\label{injected-real1}
	p_m^t = p_{r,m}^t - p_{L,m}^t, \\
	q_m^t = q_{r,m}^t - q_{L,m}^t.\label{injected-imag1}
	\end{eqnarray}
	For any node $ m \notin \mathcal{D} $, we have that $p_{r,m}^t = 0$ and $q_{r,m}^t= 0$. For notational simplicity, we assume that 
	one DER is connected at a node of the power system, and we do not describe additional constraints governing the operation of energy storage systems and electric vehicles; however, the methodology proposed in the following sections is straightforwardly applicable to the case where multiple DERs are connected at a node and where the OWPF problem includes the constraints for, e.g., the state of charges of energy storage systems  and electric vehicles.
	
	
	\section{The OWPF Problem}
	
	Water and power networks are coupled because the power consumed by a pump is proportional to the pump's pressure gain times its flow rate. Assume that a pump $ij \in \P$  is electrically connected to one node in the power network. Let $\sigma: \mathcal{P} \rightarrow \mathcal{M}$ map water pumps to the electrical nodes to which they are connected; for example, if pump $ij \in \mathcal{P}$ is connected to node $m \in \mathcal{M}$ of the distribution system, then $m = \sigma(ij)$.  For the electrical nodes $\sigma(\mathcal{P})$, we can substitute~\eqref{eq:pump-power-consumption} into~\eqref{injected-real1} to obtain:
	\begin{eqnarray}\label{injected-real1_pump}
	p_{\sigma(ij)}^t = p_{r,\sigma(ij)}^t - p_{L,\sigma(ij)}^t -   \frac{\rho g}{\eta_{ij}} \hat{h}_{ij}^t q_{ij}^t  
	\end{eqnarray}
	which couples optimization variables in both the water and power networks.
	
	With~\eqref{injected-real1_pump} in place, we formulate OWPF as the following multi-period optimization problem:
	\begin{subequations}\label{eq:OWPF}
		\begin{align}
		& \quad \min \quad \sum_{t=1}^{T} \left( C_{s}^t({\bf v}^t, p_0^t) + \sum_{m \in \mathcal{D}} C_{c,m}^t(p_{r,m}^t , q_{r,m}^t)  + \sum_{ij \in \P} C_{w,\sigma(ij)}^t (p_{\sigma(ij)}^t)   \right)
		\label{eq:OWPF:obj}
		\\ 
		&\quad \text{over}\quad \{{\bf v}^t, \{p_{r,m}^t , q_{r,m}^t\}_{\forall m \in \mathcal{D}}\}_{\forall t = 1, \ldots, T}  \nonumber
		\\
		&\quad ~~~~~\quad \{q^t, h^t, \hat h^t, \tilde h^t, w^t, R^t\}_{\forall t = 1, \ldots, T}  \nonumber
		\\
		&\quad  \text{~s.t.}  ~\quad
		\text{\eqref{eq:junction-conservation}--\eqref{eq:Darcy-Weisbach}}, 
		\text{\eqref{power-real}--\eqref{power-imag}, \eqref{setpoints}, \eqref{injected-imag1}}  \nonumber
		\\
		&  ~~~~~~~~~~   p_{m}^t = p_{r,m}^t - p_{L,m}^t -   \frac{\rho g}{\eta_{ij}} \hat{h}_{ij}^t q_{ij}^t\ \ ~\forall\! \, m\!\in\! \sigma(\mathcal{P})  \label{eq:OWPF:powerbalancePumps} \\
		&  ~~~~~~~~~~~ p_{m}^t = p_{r,m}^t - p_{L,m}^t , ~~~~~~~ \forall \, m \in \mathcal{M}^+ \backslash \sigma(\mathcal{P})  \label{eq:OWPF:powerbalanceNoPupms} \\
		&\quad  ~~~~~  \quad
		v_0^t = v_0^\text{ref},\quad \forall t = 1,\dots,T \label{eq:OWPF:substation}
		\\
		&\quad  ~~~~~  \quad
		\underline{ v}^2   \leq {{\bf v}^t}^H {\bf M}_m {\bf v}^t \leq \overline{v}^2, \quad \forall m \in \M^+, ~\forall t
		\label{eq:OWPF:voltage} 
		\end{align}
	\end{subequations}
	where the constraint~\eqref{eq:OWPF:substation} specifies a constant voltage magnitude $v_0^\text{ref}$ at the power substation, and~\eqref{eq:OWPF:voltage} imposes the voltage regulation requirement for the other nodes (i.e., ANSI C84.1 limits). Additional engineering constraints for the power network such as ampacity limits or flow limits can be added without requiring modifications to the procedure outlined in the next section~\cite{Zamzam-2017}.
	
	The cost function~\eqref{eq:OWPF:obj} is composed of three terms:
	\begin{enumerate}
		\item {\it Power system operator cost} ($ C_s^t $): this function captures operational objectives of the power network operator; for example, minimization of voltage deviations, minimization of power losses, or deviations from given setpoints for the power at the substation. For example, the latter can be expressed as $ c_s (p_0^t - p_{0,\textnormal{set}}^t)^2$, where $ c_s > 0$ is a given constant/price. 
		\item {\it DER cost} ($ C_{c,m}^t $): this function captures DER-related objectives. For example, for a PV system, one might consider the cost $\beta^t (\overline{p}_m^t - p_{r,m}^t) + R_c(p_{r,m}^t) $, where \revs{$ \overline{p}_m^t $ denotes the active power available from a PV system located at node $ m $ at time $ t $ (based on prevailing irradiance conditions)}, $ \beta^t $ is the  price of power or reward for ancillary service provisioning (received by the customers from the utility) at time $ t $, and the function $R_c(p_{r,m}^t) $ represents a convex  regularization term that can be used to promote solutions with specific characteristics.
		\item {\it Water network cost} ($ C_{w,\sigma(ij)}^t $): this function models objectives of the water network operator.  For example, payments of the water network operator for the power consumed by the pumps can be expressed as $ \alpha^t p_{\sigma(ij)}^t + R_w(p_{\sigma(ij)}^t) $, where $\alpha^t$ is the price of electricity (received from the power systems operator) and the regularization term $ R_w(p_{\sigma(ij)}^t) $ is used again to promote the preferred features of the optimal solution (e.g., smoothness or sparsity). In particular, the sparsity of the solution leads to the case where only a subset of the pumps are used. On the other hand, the smoothness of the solution prevent the case where one of the pumps is overloaded while the other pumps do not operate.
	\end{enumerate}
	In the OWPF problem~\eqref{eq:OWPF}, the coupling across multiple time periods appears only in the tank dynamics \eqref{eq:tank-dynamics} (with  $h_k^{0,\textnormal{out}}$ in \eqref{eq:tank-dynamics} a given constant); however, additional constraints for the  state of charge of  energy storage systems  and electric vehicles can be straightforwardly added to the problem formulation. 
	
	Problem~\eqref{eq:OWPF} is \emph{nonconvex} because of the constraints~\eqref{eq:pump-head-gain}, \eqref{eq:Darcy-Weisbach}, \eqref{power-real}, \eqref{power-imag},~\eqref{eq:OWPF:powerbalancePumps}, and the lower bound in~\eqref{eq:OWPF:voltage}; however,  these constraints will be reformulated as  nonconvex  quadratic inequalities in the following section, and they will  be efficiently managed by the FPP-SCA {approach}~\cite{Mehanna-2015,Zamzam-2017}.

	\section{Successive Convex Approximation}
	
	The FPP-SCA is a two-step algorithm that solves convex problems iteratively. In the first step, a feasible solution is obtained by solving  an inner approximation of~\eqref{eq:OWPF}. The inner approximation is obtained by rewriting a nonconvex quadratic inequality as a difference of two convex functions, and then linearizing the concave term around a given restriction point.  To promote the feasibility of the approximation, we add slack variables  to the constraints and minimize the norm of the slack variables over the approximate feasible set. We then use the solution as a restriction point for the next step. If the slack variables become zero, we by construction have a feasible point of~\eqref{eq:OWPF} . In the second stage, we solve a sequence of convex inner approximations
	of~\eqref{eq:OWPF} until convergence to a KKT point.

	
	\subsection{Nonconvexity in Power Network Constraints}
	Notice first that the matrices $ {\bf Y}_m $ and $ \overline{\bf Y}_m $ are indefinite~\cite{low2014convex}. Rewrite~\eqref{power-real} as the following two inequalities:
	\begin{eqnarray}
	{{\bf v}^t}^H {\bf Y}_m {\bf v}^t \leq p_m^t,\label{power-real1}\\
	{{\bf v}^t}^H ({-\bf Y}_m) {\bf v}^t \leq -p_m^t.\label{power-real2}
	\end{eqnarray}
	Constraint~\eqref{power-real1} can be rewritten as: 
	\begin{equation}
	{{\bf v}^t}^H ({\bf Y}_m^{(+)}+ {\bf Y}_m^{(-)}) {\bf v}^t \leq p_m^t
	\end{equation}
	where $ {\bf Y}_m^{(+)} $ and $ {\bf Y}_m^{(-)} $ are the positive and negative semidefinite parts of matrix $ {\bf Y}_m $ {(obtained through eigenvalue decomposition)}. Focusing on the negative semidefinite matrix $ {\bf Y}_m^{(-)} $, we can write: 
	\begin{equation}
	({\bf v}^t - {\bf u}^t)^H {\bf Y}_m^{(-)} ({\bf v}^t - {\bf u}^t) \leq 0
	\end{equation}
	where $ {\bf u}^t $ is a  {\em restriction point}. By rearranging terms, we get:
	\begin{equation}
	\label{power-real_der}
	{{\bf v}^t}^H {\bf Y}_m^{(-)} {\bf v}^t \leq 2 \Re\{{{\bf u}^t}^H {\bf Y}_m^{(-)} {\bf v}^t\} - {{\bf u}^t}^H {\bf Y}_m^{(-)} {\bf u}^t.
	\end{equation}
	Using~\eqref{power-real_der}, a convex surrogate of the nonconvex constraint~\eqref{power-real1} can be obtained as follows:
	\begin{align}
	& {{\bf v}^t}^H {\bf Y}_m^{(+)} {\bf v}^t + 2 \Re\{{{\bf u}^t}^H {\bf Y}_m^{(-)} {\bf v}^t\} \leq p_m^t + {{\bf u}^t}^H {\bf Y}_m^{(-)} {\bf u}^t + s_{m,t}^{(p)} \label{power-real-surr1}
	\end{align}
	where the non-negative slack variable $ s_{m,t}^{(p)} $ is utilized to ensure the feasibility of the constraint. Similarly, the nonconvex constraint~\eqref{power-real2} is replaced by:
	\begin{align}
	& - {{\bf v}^t}^H {\bf Y}_m^{(-)} {\bf v}^t - 2 \Re\{{{\bf u}^t}^H {\bf Y}_m^{(+)} {\bf v}^t\}\leq -p_m^t - {{\bf u}^t}^H {\bf Y}_m^{(+)} {\bf u}^t + s_{m,t}^{(p)}. \label{power-real-surr2}
	\end{align}
	
	Following a similar procedure, the quadratic equality constraint~\eqref{power-imag} can be replaced by the following convex constraints:
	\begin{align}
	&	{{\bf v}^t}^H \overline{\bf Y}_m^{(+)} {\bf v}^t + 2 \Re\{{{\bf u}^t}^H \overline{\bf Y}_m^{(-)} {\bf v}^t\} \leq q_m^t + {{\bf u}^t}^H \overline{\bf Y}_m^{(-)} {\bf u}^t + s_{m,t}^{(q)}, \label{power-imag-surr1}
	\end{align}
	\begin{align}\
	&		- {{\bf v}^t}^H {\bf \overline{Y}}_m^{(-)} {\bf v}^t - 2 \Re\{{{\bf u}^t}^H \overline{\bf Y}_m^{(+)} {\bf v}^t\}\leq -q_m^t - {{\bf  u}^t}^H \overline{\bf Y}_m^{(+)} {\bf u}^t + s_{m,t}^{(q)} \label{power-imag-surr2}
	\end{align} 
	where the $ s_{m,t}^{(q)} $ is  added in order to ensure feasibility of the convex inner-approximation, while ${\bf \overline{Y}}_m^{(+)}$ and ${\bf \overline{Y}}_m^{(-)}$ are the positive and negative semidefinite parts of  $ {\bf \overline{Y}}_m $.
	
	Finally, the lower bound on the bus voltage magnitude~\eqref{eq:OWPF:voltage} is replaced by: 
	\begin{equation}\label{power-vol-surr}
	- 2 \Re\{{{\bf u}^t}^H {\bf M}_m {\bf v}^t\} \leq -\underline{v}^2 - {{\bf u}^t}^H {\bf M}_m {\bf u}^t + s_{m,t}^{(v)},
	\end{equation}
	where $ s_{m,t}^{(v)}$ is a non-negative slack variable.

	\subsection{Nonconvexity in Water Network Constraints}
	Constraint~\eqref{eq:pump-head-gain} is nonconvex. However, it has been shown in~\cite{fooladivanda2017energy}  that~\eqref{eq:pump-head-gain}  can be  replaced by the following constraint without loss of optimality: 
	\begin{IEEEeqnarray}{rCl}
		\hat{h}^t_{ij} &\leq&  A_{ij} \left(q_{ij}^t\right)^2 + B_{ij}^{\text{max}} q_{ij}^t  + C_{ij}^{\text{max}} 
		,  \label{eq:pump-head-gain-reformulated}
	\end{IEEEeqnarray}
	where $B_{ij}^{\text{max}}:= B_{ij} w_{ij}^{\text{max}}$ and $C_{ij}^{\text{max}}  := C_{ij} \left(w_{ij}^\text{max}\right)^2$. 
	Because $A_{ij} < 0$,  \eqref{eq:pump-head-gain-reformulated} is a convex quadratic constraint. With this reformulation, the variables $w_{ij}^t$ can be eliminated; in fact, a unique $w_{ij}^t$ that satisfies \eqref{eq:max-speed} can always be recovered from any feasible solution $(q_{ij}^t, \hat h_{ij}^t)$. 
	
	Next, consider replacing~\eqref{eq:Darcy-Weisbach} with the following two inequalities:
	\begin{IEEEeqnarray}{rCl} \label{eq:Darcy-Weisbach-reformulated-1}
		\tilde{h}^t_{ij} \!&\geq&  F_{ij} \left(q_{ij}^t\right)^2 \label{eq:Darcy-Weisbach-reformulated-11} \\
		-\tilde{h}^t_{ij} \!&\geq& \! -F_{ij} \left(q_{ij}^t\right)^2 \label{eq:Darcy-Weisbach-reformulated-12}
	\end{IEEEeqnarray}
	and notice that~\eqref{eq:Darcy-Weisbach-reformulated-11} is convex because $F_{ij} > 0$~\cite{fooladivanda2017energy}. On the other hand,~\eqref{eq:Darcy-Weisbach-reformulated-12} is nonconvex. A convex restriction of~\eqref{eq:Darcy-Weisbach-reformulated-12} can be written as:
	\begin{equation}\label{eq:Darcy-Weisbach-reformulated-2}
	F_{ij} (z_{q,ij}^t)^2 - 2z_{q,ij}^t F_{ij} q_{ij}^t \leq -\tilde{h}_{ij}^t + s_{ij,t}^{(q)}
	\end{equation}
	where the non-negative slack variable $ s_{ij,t}^{(q)} $ is added to insure feasibility.
	
	Regarding the power balance constraint~\eqref{eq:OWPF:powerbalancePumps},  it is convenient to introduce the auxiliary variables  $ \beta_{ij}^t $ and $p_{\textrm{pump},ij}^t $ per pump $ij \in \mathcal{P}$, and rewrite~\eqref{eq:OWPF:powerbalancePumps} as the following equivalent set of constraints:
	\begin{align}
	p_{\sigma(ij)}^t & = p_{r,\sigma(ij)}^t - p_{L,\sigma(ij)}^t - p_{\textrm{pump},\sigma(ij)}^t  \label{eq:OWPF:powerbalancePumps_2} \\ 
	\beta_{ij}^t & = p_{\textrm{pump},ij}^t \label{eq:eq:OWPF:powerbalancePumps_aux} \\
	\beta_{ij}^t & = \frac{\rho g}{\eta_{ij}} \hat{h}_{ij}^t q_{ij}^t . \label{eq:OWPF:powerbalancePumps_aux2}
	\end{align}
	The auxiliary variables  $ \beta_{ij}^t $ and $p_{\textrm{pump},ij}^t $ facilitate the application of the FPP-SCA method as well as the development of a distributed algorithm in Section~\ref{sec:Distributed}.
	Re-express~\eqref{eq:OWPF:powerbalancePumps_aux2} as:
	\begin{IEEEeqnarray}{rCl}
		{{\bf y}_{ij}^t}^T {\bf T}_{ij} {\bf y}_{ij}^t \leq - \beta_{ij}^t   \label{eq:power-coupling-water3}\\
		{{\bf y}_{ij}^t}^T (-{\bf T}_{ij}) {\bf y}_{ij}^t \leq  \beta_{ij}^t  \ \label{eq:power-coupling-water4}
	\end{IEEEeqnarray}
	where ${\bf y}_{ij}^t := [\hat{h}_{ij}^t, q_{ij}^t]^T$ for simplicity and $ {\bf T}_{ij} $ is a two-by-two matrix with  $ \frac{\rho g}{2\eta_{ij}} $ on the off-diagonal entries. Constraints~\eqref{eq:power-coupling-water3}--\eqref{eq:power-coupling-water4} are nonconvex. Rewriting $ {\bf T}_{ij}$ as $ {\bf T}_{ij} =  {\bf T}_{ij}^{(+)} +  {\bf T}_{ij}^{(-)}$, where  ${\bf T}_{ij}^{(+)}$ and \revs{${\bf T}_{ij}^{(-)}$} are the positive and negative semidefinite {parts of  ${\bf T}_{ij}$, surrogate} convex constraints for~\eqref{eq:power-coupling-water3} and~\eqref{eq:power-coupling-water4} can be written as: 
	\begin{align}
	&	{{\bf y}_{ij}^t}^T {\bf T}_{ij}^{(+)} {\bf y}_{ij}^t + 2 {{\bf z}_{ij}^{t}}^T {\bf T}_{ij}^{(-)} {\bf y}_{ij}^{t}\leq - \beta_{ij}^t + {{\bf z}_{ij}^{t}}^T {\bf T}_{ij}^{(-)} {\bf z}_{ij}^{t} + s_{ij,t}^{(p)}, \label{eq:power-coupling-water5}
	\end{align}
	\begin{align}
	&	-{{\bf y}_{ij}^t}^T {\bf T}_{ij}^{(-)} {\bf y}_{ij}^t - 2 {{\bf z}_{ij}^{t}}^T {\bf T}_{ij}^{(+)} {\bf y}_{ij}^{t}\leq \beta_{ij}^t - {{\bf z}_{ij}^{t}}^T {\bf T}_{ij}^{(+)} {\bf z}_{ij}^{t} + s_{ij,t}^{(p)}    \label{eq:power-coupling-water6}
	\end{align}
	\revs{where $ {\bf z}_{ij}^{t} $ represents any linearization (restriction) point, and} $s_{ij,t}^{(p)}$ is a non-negative slack variable. 
	
	We are now ready to outline the FPP-SCA algorithm for solving OWPF -- the subject of the next section. 
	
	\subsection{FPP-SCA Algorithm}
	
	For notational simplicity, collect  in the matrices $ {\bf S}_p, {\bf S}_q, {\bf S}_v \in \mathbb{R}_+^{M\times T}$   the slack variables $ s_{m,t}^{(p)} $, $ s_{m,t}^{(q)} $, and $ s_{m,t}^{(v)} $, respectively. Similarly, let the matrices $ \overline{\bf S}_p \in \mathbb{R}_+^{|\P|\times T}$ and $ \overline{\bf S}_q \in \mathbb{R}_+^{n\times T}$, where $ n = |\L \textfractionsolidus (\P \cup \V)| $, collect all the  slack variables $ s_{ij,t}^{(p)} $ and $ s_{ij,t}^{(q)} $, respectively. In addition, let $ {\bf Q}^t \in  \mathbb{R}_+^{|\L|}$, $ {\boldsymbol \beta}^t \in \mathbb{R}^{|\P|} $,  and $ {\bf h}^t \in  \mathbb{R}^{|\N|} $ collect  all the water flow variables, powers consumed by pumps, and head pressures. Lastly,  let the matrices $ {\bf Z} \in \mathbb{R}^{2|\P|\times T}$, ${\bf Z}_q \in \mathbb{R}^{n\times T}$, and ${\bf U}\in \mathbb{C}^{M\times T} $ collect the restriction points $ {\bf z}_{ij}^t $, $ z_{q,ij}^t $, and $ {\bf u}^t $, respectively. 
	
	With this notation in place, define the convex sets $ \Psi_{\bf U} $ and $ \Omega_{{\bf Z}, {\bf Z}_q} $ (which are functions of the restriction points) pertaining to power and water networks, respectively, as: 
	\begin{multline}\label{eq:FPP-power-set}
	\Psi_{\bf U} := \left\{ (\{ {\bf v}^t, \bp_{r}^t, \bq_{r}^t, \bp_{\textrm{pump}}^t  \}_{t=1}^T, {\bf S}_p, {\bf S}_q, {\bf S}_v) \ |
	\begin{split}
	&\quad\quad\quad \text{\eqref{setpoints}, \eqref{injected-imag1},~\eqref{eq:OWPF:powerbalanceNoPupms}, \eqref{eq:OWPF:substation}, \eqref{power-real-surr1}--\eqref{power-vol-surr},~\eqref{eq:OWPF:powerbalancePumps_2} } \\
	& {{\bf v}^t}^H {\bf M}_m {\bf v}^t \leq \overline{v}^2 \quad \forall m\in \M^+, \forall t = 1, \ldots, T
	\end{split}
	\right\},
	\end{multline}
	\begin{multline}\label{eq:FPP-water-set}
	\Omega_{{\bf Z}, {\bf Z}_q} :=  \bigg\{ (\{ {\bf Q}^t, {\bf h}^t, \boldsymbol{\beta}^t \}_{t=1}^T, \overline{\bf S}_p, \overline{\bf S}_q)\ |\ 
	\text{\eqref{eq:junction-conservation}--\eqref{eq:zero-head-reservoir}, \eqref{eq:on-pump-head-constraint}--\eqref{eq:on-pipe-flow-rate}},
	\text{ \eqref{eq:Darcy-Weisbach-reformulated-1}, \eqref{eq:Darcy-Weisbach-reformulated-2}, \eqref{eq:power-coupling-water5}, \eqref{eq:power-coupling-water6} }
	\bigg\}
	\end{multline}
	where $\bp_{r}^t := [\{p_{r,m}^t\}_{m \in \mathcal{D}}]^T$, $\bq_{r}^t := [\{q_{r,m}^t\}_{m \in \mathcal{D}}]^T$, and $\bp_{\textrm{pump}}^t  := [\{p_{\textrm{pump},m}^t \}_{m \in \sigma(\mathcal{P})}]^T$. \\

	It follows that the convex optimization problem to be solved at the $ k $-th iteration of the feasibility phase of the FPP-SCA algorithm can be compactly written as:
	\begin{subequations}\label{eq:OWPF-f}
		\begin{align}
		& \quad \min \quad  \| {\bf S}_p\|_F^2 + \| {\bf S}_q\|_F^2 + \| {\bf S}_v\|_F^2 + \| \overline{\bf S}_p\|_F^2 + \| \overline{\bf S}_q\|_F^2
		\label{eq:OWPF-f:obj}
		\\ 
		&\quad \text{over}\quad \{ {\bf v}^t, \bp_{r}^t, \bq_{r}^t, \bp_{\textrm{pump}}^t \}_{t=1}^T, {\bf S}_p, {\bf S}_q, {\bf S}_v  \nonumber
		\\
		&\quad ~~~~~\quad \{ {\bf Q}^t, {\bf h}^t, \boldsymbol{\beta}^t \}_{t=1}^T, \overline{\bf S}_p, \overline{\bf S}_q \nonumber
		\\
		&\quad  \text{~s.t.}  ~\quad
		(\{ {\bf v}^t, \bp_{r}^t, \bq_{r}^t, \bp_{\textrm{pump}}^t \}_{t=1}^T, {\bf S}_p, {\bf S}_q, {\bf S}_v) \in \Psi_{{\bf U}{(k)}}\\
		&\quad  ~~~~~  \quad
		(\{ {\bf Q}^t, {\bf h}^t, \boldsymbol{\beta}^t \}_{t=1}^T, \overline{\bf S}_p, \overline{\bf S}_q) \in \Omega_{{\bf Z}{(k)}, {\bf Z}_q{(k)}}  \nonumber
		\\
		&\quad  ~~~~~  \quad
		\beta_{ij}^t  = p_{\textrm{pump},\sigma(ij)}^t  \quad \forall ij \in \P, t = 1, \ldots, T 
		\end{align}
	\end{subequations}
	where the cost function aims to minimize the violation of the nonconvex constraint in the original OWPF, and ${\bf U}(k)$, ${\bf Z}(k), {\bf Z}_q(k)$ are the restriction points at iteration $k$. In particular, ${\bf U}(k)$, ${\bf Z}(k), {\bf Z}_q(k)$ coincide with the optimal solution of~\eqref{eq:OWPF-f} at the previous iteration $k-1$. Algorithm~\ref{Algorithm: FPP-Alg} summarizes the overall feasibility procedure. 
	Notice that problem~\eqref{eq:OWPF-f} is a second-order cone program, and observe that the cost function is monotonically nonincreasing because the restriction points are always feasible. Although
	this method in not guaranteed to find a feasible point, it converged in all simulations to operating points that were feasible for the water and power systems {considered}.

	\SetKw{Init}{Initialization: }
	\DontPrintSemicolon
	\begin{algorithm}[h]
		\renewcommand{\arraystretch}{1.2}
		\caption{FPP Method for OWPF}
		{\footnotesize
			\Init{\rm set $ k = 0 $ and $ {\bf U}{(0)} $ to be all ones, and $ {\bf Z}{(0)}, {\bf Z}_q{(0)} $ to be all zeros.}\;
			\Repeat{$(\| {\bf S}_p\|_F^2 + \| {\bf S}_q\|_F^2 + \| {\bf S}_v\|_F^2 + \| \overline{\bf S}_p\|_F^2 + \| \overline{\bf S}_q\|_F^2)  \leq \epsilon $}{
				
				[S1]	Solve~\eqref{eq:OWPF-f}: 
				$
				\{\{ {\bf v}^t(k), \bp_{r}^t(k), \bq_{r}^t(k), \bp_{\textrm{pump}}^t(k) \}_{t=1}^T,
				 \{ {\bf Q}^t(k), {\bf h}^t(k), \boldsymbol{\beta}^t (k)\}_{t=1}^T,
				{\bf S}_p(k), {\bf S}_q(k), {\bf S}_v(k), \overline{\bf S}_p(k), \overline{\bf S}_q(k)  \leftarrow
				\text{solution of \eqref{eq:OWPF-f}} \nonumber
				$\;

				[S2]	Update restriction points: 
				
				$  {\bf U}{(k+1)}  \leftarrow  \{{\bf v}^t(k)\}_{t=1}^T $\;
				$ {\bf Z}{(k+1)} \leftarrow  \{q_{ij}^t(k), \hat{h}_{ij}^t(k), ij \in \mathcal{P}\}_{t=1}^T $\;
				$ {\bf Z}_q{(k+1)}  \leftarrow \{q_{ij}^t(k),  ij \in \L \textfractionsolidus (\P \cup \V)\}_{t=1}^T $ 
				
				\vspace{.2cm}		
				
				[S3]	Increase iteration index:
				
				$ k \leftarrow k+1 $.\;
				
				\vspace{.2cm}	
				
			}
		}
		\label{Algorithm: FPP-Alg}
	\end{algorithm}
	
	Once a feasible operating point is identified through Algorithm~\ref{Algorithm: FPP-Alg}, the nonconvex feasible set of OWPF is replaced at each iteration by an inner convex approximation. The optimization problem to be solved at the $ k $-th iteration of the refinement phase can be formulated as:
	\begin{subequations}\label{eq:OWPF-r}
		\begin{align}
		& \quad \min \quad \sum_{t=1}^{T} \left( C_{s}^t({\bf v}^t, p_0^t) + \sum_{m \in \mathcal{D}} C_{c,m}^t(p_{r,m}^t , q_{r,m}^t) + \sum_{ij \in \P} C_{w,\sigma(ij)}^t (p_{\sigma(ij)}^t)   \right)
		\label{eq:OWPF-r:obj}
		\\ 
		&\quad \text{over}\quad \{ {\bf v}^t, \bp_{r}^t, \bq_{r}^t, \bp_{\textrm{pump}}^t,  {\bf Q}^t, {\bf h}^t, \boldsymbol{\beta}^t \}_{t=1}^T \nonumber
		\\
		&\quad  \text{~s.t.}  ~\quad
		(\{ {\bf v}^t, \bp_{r}^t, \bq_{r}^t, \bp_{\textrm{pump}}^t \}_{t=1}^T) \in \tilde{\Psi}_{{\bf U}{(k)}}\\
		&\quad  ~~~~~  \quad
		(\{ {\bf Q}^t, {\bf h}^t, \boldsymbol{\beta}^t \}_{t=1}^T) \in \tilde{\Omega}_{{\bf Z}{(k)}, {\bf Z}_q{(k)}}  \nonumber
		\\
		&\quad  ~~~~~  \quad
		\label{eq:OWPF-r:coupling}
		\beta_{ij}^t  = p_{\textrm{pump},\sigma(ij)}^t  \quad \forall ij \in \P, t = 1, \ldots, T 
		\end{align}
	\end{subequations}
	where $ \tilde{\Psi}_{{\bf U}{(k)}} $ is the restriction of the set $ {\Psi}_{{\bf U}{(k)}} $ to the plane $ {\bf S}_p = {\bf S}_q = {\bf S}_v = {\boldsymbol 0} $. Similarly, the set $ \tilde{\Omega}_{{\bf Z}{(k)}, {\bf Z}_q{(k)}} $ is the projection of $ {\Omega}_{{\bf Z}{(k)}, {\bf Z}_q{(k)}} $ when $ \overline{\bf S}_q = \overline{\bf S}_p = {\boldsymbol 0} $. At each step $k = 2, \ldots$, the sets $ \tilde{\Psi}_{{\bf U}{(k)}} $ and $ \tilde{\Omega}_{{\bf Z}{(k)}, {\bf Z}_q{(k)}} $ are formed based on the optimal solution of~\eqref{eq:OWPF-r} at the previous iteration $k - 1$.  The proposed methodology generates a monotone sequence that converges
	to a KKT point of the original OWPF~\eqref{eq:OWPF}, and it is given as Algorithm~\ref{Algorithm: SCA-Alg}.

	\SetKw{Init}{Initialization: }
	\DontPrintSemicolon
	\begin{algorithm}[h]
		\renewcommand{\arraystretch}{1.2}
		\caption{SCA Method for OWPF}
		{\footnotesize
			\Init{\rm set $ k = 0 $ and construct $ {\bf U}{(0)}, {\bf Z}{(0)},  {\bf Z}_q{(0)} $ based on the solution of Algorithm~\ref{Algorithm: FPP-Alg}.}\;
			\Repeat{Cost reduction $\leq \epsilon $}{
				
				[S1]	Solve~\eqref{eq:OWPF-r}: 	
				
				$
				\{\{ {\bf v}^t(k), \bp_{r}^t(k), \bq_{r}^t(k), \bp_{\textrm{pump}}^t(k) \}_{t=1}^T,
				\{ {\bf Q}^t(k), {\bf h}^t(k), \boldsymbol{\beta}^t (k)\}_{t=1}^T \leftarrow
				\text{solution of \eqref{eq:OWPF-r}} \nonumber
				$\;
				
				[S2]	Update restriction points: 
				
				$  {\bf U}{(k+1)}  \leftarrow  \{{\bf v}^t(k)\}_{t=1}^T $\;
				$ {\bf Z}{(k+1)} \leftarrow  \{q_{ij}^t(k), \hat{h}_{ij}^t(k), ij \in \mathcal{P}\}_{t=1}^T $\;
				$ {\bf Z}_q{(k+1)}  \leftarrow \{q_{ij}^t(k), ij \in \L \textfractionsolidus (\P \cup \V)\}_{t=1}^T $ \;
				
				\vspace{.2cm}		
				
				[S3]	Increase iteration index:
				
				$ k \leftarrow k+1 $.\;
				
				\vspace{.2cm}	
				
			}
		}
		\label{Algorithm: SCA-Alg}
	\end{algorithm}

	Algorithm~\ref{Algorithm: SCA-Alg} has converged if the cost function between two consecutive steps is less than a given quantity $0< \epsilon \ll 1$.
	\revs{
		\begin{claim}[Convergence] From~\cite[Theorem~1]{razaviyayn-2014}, it follows that every limit point generated by the proposed algorithms is a KKT point. Hence, the first phase converges to a KKT point of \eqref{eq:OWPF-f}. In addition, if the algorithm starts  the second phase from a feasible point, then the sequence  converges to the set containing all the KKT points of the problem \eqref{eq:OWPF}.
		\end{claim}
	}
	\revs{ The convergence of the first phase of the algorithm follows from~\cite{razaviyayn-2014}. 
		Since the SCA phase is initialized from a feasible point of the OWPF, the sequence produced by Algorithm~\ref{Algorithm: SCA-Alg} is always feasible; i.e., every point in the sequence lies in the feasibility set of~\eqref{eq:OWPF}. Thus, the  sequence converges to a set that contains all the KKT points of~\eqref{eq:OWPF}.
	}

	In the next section, a distributed solution strategy will be outlined  based on the ADMM.

	\section{Distributed Algorithm}
	\label{sec:Distributed}
	Power and water networks are usually planned and operated independently. In this section, we develop a distributed solver where  water and power operators pursue individual operational objectives and retain controllability of their own assets (pumps for the water network and DERs for the power network).   In particular,  in the proposed setting, each system operator solves a smaller optimization problem with variables associated only with its network and controllable devices -- a rendition of the AC optimal power flow  for the power network operator, and an optimal water flow problem for the water system operator. The two operators subsequently  exchange information regarding the shared variables to \emph{reach consensus} on the powers consumed by the water pumps. Using this approach, both systems keep their private information and only exchange the optimized values of the powers consumed by the pumps. 
	
	Assume that each system has determined an operational profile that satisfies all network-level constraints. The proposed technical approach consists of defining the following augmented Lagrangian associated with~\eqref{eq:OWPF-r}: 
	\begin{multline}\label{eq:OWPF-Lagrange}
	L_\rho (\left\{ {\bf v}^t, \bp_{r}^t, \bq_{r}^t , \bp_{\textrm{pump}}^t, {\bf Q}^t, {\bf h}^t, \boldsymbol{\beta}^t,  \{\nu_{ij}^t\}_{ij\in \P} \right\}_{t=1}^T) = \sum_{t=1}^{T} \left( C_{s}^t(p_0^t) + \sum_{m \in \mathcal{D}} C_{c,m}^t(p_m^t) + \sum_{ij \in \P} C_{w,\sigma(ij)}^t(\beta_{ij}^t)\right. \\ 
	\left. +\sum_{ij\in\mathcal{P}} \nu_{ij}^t \left(\beta_{ij}^t - p_{\textrm{pump},\sigma(ij)}^t \right) + \frac{\rho}{2} \left(\beta_{ij}^t - p_{\textrm{pump},\sigma(ij)}^t \right)^2 \right)
	\end{multline}
	where $ \nu_{ij}^t $ is the dual variable associated with the consensus-enforcing constraint~\eqref{eq:OWPF-r:coupling} and $ \rho > 0$ is an given  parameter. 
	
	Leveraging the decomposability of \eqref{eq:OWPF-Lagrange}, the proposed ADMM-based algorithm is summarized in Algorithm~\ref{Algorithm: Dist-ADMM}, wherein $k $ is the iteration index.
	
	\SetKw{Init}{Initialization: }
	\DontPrintSemicolon
	\begin{algorithm}
		{\footnotesize
			\Init{set $ k = 0 $}, $ {\nu_{ij}^t}^{(0)} = {0}\ \forall ij\in\P, t = 1, \ldots, T$, and each system initializes its own variable from a feasible operational profile.\;
			\Repeat{convergence criterion is met}{
				{
					[S1a] Update power-related variables: 
					$$
					(\{ {\bf v}^t(k+1), \bp_{r}^t(k+1), \bq_{r}^t(k+1), \bp_{\textrm{pump}}^t(k+1)  \}_{t=1}^T)) \leftarrow\ 
					\underset{\{ {\bf v}^t, \bp_{r}^t,  \bq_{r}^t , \bp_{\textrm{pump}}^t\} \in \tilde{\Psi}_{{\bf U}{(k)}}}{\arg\min} L_\rho (\left\{ {\bf v}^t, \bp_{r}^t, \bq_{r}^t, \bp_{\textrm{pump}}^t,  \{\nu_{ij}^t(k)\} \right\}_{t=1}^T)  \nonumber 
					$$   
				}
				\vspace{0.1cm}
				
				[S1b] Update water-related variables: 		
				$$
				(\{{\bf Q}^t(k+1), {\bf h}^t(k+1), \boldsymbol{\beta}^t(k+1) \}_{t=1}^T) \leftarrow
				\underset{(\{ {\bf Q}^t, {\bf h}^t, \boldsymbol{\beta}^t \}_{t=1}^T) \in \tilde{\Omega}_{{\bf Z}{(k)}, {\bf Z}_q{(k)}}}{\arg\min} L_\rho (\left\{\{\nu_{ij}^t(k)\}, {\bf Q}^t,\! {\bf h}^t,\! \boldsymbol{\beta}^t\! \right\}_{t=1}^T) \nonumber
				$$ 
				
				[S2] Update dual  variables: 	
				
				\For{
					$ ij \in \P, t = 1, \ldots, T$
				}
				{
					$ {{\nu}_{ij}^t}{(k+1)} = {{\nu}_{ij}^t}{(k)} + \rho \left({{\beta}_{ij}^t}{(k+1)} - {p_{\textrm{pump},\sigma(ij)}^t}{(k+1)}\right)  $ \;
				}
				\vspace{0.1cm}
				
				[S3] Update restriction points:
				
				$  {\bf U}{(k+1)}  \leftarrow  \{{\bf v}^t(k+1)\}_{t=1}^T $\;
				$ {\bf Z}{(k+1)} \leftarrow  \{q_{ij}^t(k+1), \hat{h}_{ij}^t(k+1), ij \in \mathcal{P}\}_{t=1}^T $\;
				$ {\bf Z}_q{(k+1)}  \leftarrow \{q_{ij}^t(k+1),  ij \in \L \textfractionsolidus (\P \cup \V)\}_{t=1}^T $
				\vspace{0.1cm}	
				
				[S4] Increase iteration index:
				$ k \leftarrow k+1 $ \;
			}
			\caption{Distributed OWPF algorithm.}}
		\label{Algorithm: Dist-ADMM}
	\end{algorithm}
	
	{
		Each iteration of the procedure described in Algorithm \ref{Algorithm: Dist-ADMM} consists of two main steps. In step [S1] power and water operators solve two local subproblems to update variables $\{{\bf v}^t(k+1), \bp_{r}^t(k+1), \bq_{r}^t(k+1), \bp_{\textrm{pump}}^t(k+1)  \}_{t=1}^T$ and $\{{\bf Q}^t(k+1), {\bf h}^t(k+1), \boldsymbol{\beta}^t(k+1) \}_{t=1}^T$, respectively. The subproblem solved by the power network operator coincides with an AC optimal power flow, with the cost augmented by the consensus-enforcing regularization term $\sum_{t = 1}^T \sum_{m \in \sigma(\mathcal{P})} ( - \nu_{m}^t p_{\textrm{pump},m}^t + \frac{\rho}{2} (\beta_{m}^t(k) - p_{\textrm{pump},m}^t)^2)$. The problem solved by the water system operator is an optimal water flow~\cite{fooladivanda2017energy}, augmented with the  regularization term $\sum_{t = 1}^T \sum_{ij \in \mathcal{P}} ( \nu_{ij}^t \beta_{ij}^t + \frac{\rho}{2} (\beta_{ij}^t - p_{\textrm{pump},\sigma(ij)}^t(k))^2)$. After solving the two problems, the water and power operators exchange the variables $\{p_{\textrm{pump},m}^t(k+1)\}$ and $\{\beta_{m}^t(k+1) \}$ and they locally perform the dual step in [S2]. 
		\revf{The algorithm terminates when the norm of the difference between two consecutive restriction points is less than $ \epsilon $. Once the algorithm  terminates,} it holds that ${\beta}_{ij}^t = {p_{\textrm{pump},\sigma(ij)}^t}$; i.e., the two systems agree on the power consumed by the pumps.

		\section{Experimental Results}
		\label{sec:experiments}
		In this section, the proposed approach is tested in a setting where  the IEEE 13-node test feeder is coupled with a municipal water distribution network~\cite{Cohen-2000}. The OWPF is solved using the proposed centralized and the distributed solvers. The MATLAB-based optimization package YALMIP~\cite{yalmip} is used along with the interior-point solver SeDuMi~\cite{sedumi}.  The centralized solver initializes the voltage with the flat voltage profile; the parameter $ \epsilon $ is set to $ 10^{-5} $.  
		
		The water network is shown in Fig.~\ref{fig:water-network} and it is modified version of the ones described in~\cite{Cohen-2000}. It consists of $ 7 $ junctions, $ 2 $ reservoirs, and a tank and it features $ 5 $ pipes, $ 2 $ pressure-reducing valves, and $ 3 $ pumps. The elevation of all junctions and the minimum allowable head pressure are given in Table \ref{table:junctions-info}. The elevation of the water reservoirs at nodes $ 8 $ and $ 9 $ are $  -2.5 $ meters and $ 5 $ meters, respectively. The water demands over the optimization time at junctions $ 3 $, $ 4 $, $ 5 $, $ 6 $, and $ 7 $ are shown in Fig.~\ref{fig:water-demand}. The tank has an area of $ 490.87 ~\text{m}^2  $ and maximum height of $ 30 $ meters. The initial height of the water in the tank is $ 10 $ meters. As an example of application, the optimization problem is solved over an interval of $ 12 $ hours (from $ 6 $ AM to $ 6 $ PM), with intervals of $ 30 $ minutes. The regularization term of the water cost function is chosen to be the Frobenius norm of the power consumption of the pumps over the optimization period with a weight of $ 10^{-3} $. 
		
		
		\begin{figure}[t]
			\centering
			\includegraphics[width=9.5cm]{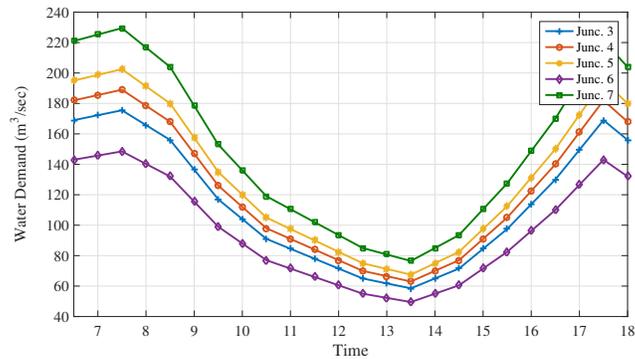}
			\caption{Water demand at the junctions with time}
			\label{fig:water-demand}
			\vspace{-.25cm}
		\end{figure}

		\begin{table}[htbp]
			\centering
			\caption{Junctions details.}
			\begin{tabular}{|c|c|c|}
				\hline
				{\textbf{Junction}} & {\textbf{Elevation} (m)} &{\textbf{Minimum Head Pressure} (m) } \\
				\hline
				1     & 6     & 0 \\
				\hline
				2     & 33    & 0 \\
				\hline
				3     & 1.5   & 35 \\
				\hline
				4     & -8    & 40 \\
				\hline
				5     & 33    & 40 \\
				\hline
				6     & 8     & 35 \\
				\hline
				7     & 4     & 40 \\
				\hline
			\end{tabular}%
			\label{table:junctions-info}%
		\end{table}%
		
		The pumping stations are connected to the IEEE 13-node distribution test feeder, which is shown in Fig. \ref{fig:power-network}. In particular, pump  $8 \rightarrow 1$ is connected to node 633 of the feeder, pump $9 \rightarrow 2$ to node 645, and pump $4 \rightarrow 5$ to node 684. A single-phase model of the distribution feeder is considered, wherein PV inverters are assumed to be located at the nodes 634, 646, 675, 611, and 652. Because of the high PV penetration, the system is likely to experience overvoltage challenges. In this case, curtailment of the active power at the PV units is necessary to maintain voltage magnitudes within  prescribed limits. The  prices of electricity (utilized in the cost function for the water network and to discourage curtailment from PV systems) are obtained from the Midcontinent Independent System Operator\footnote{Available at: \url{https://www.misoenergy.org}} for June, 20, 2017. The regularization term of the curtailment cost function is chosen to be a weighted Frobenius norm of the active power curtailed during the optimization period with a weight of $ 10^{-3} $.
		

		The OWPF strategy is compared to the decoupled case, where: 1) water system operator solves an optimal flow problem to minimize the power consumption based on the electricity prices~\cite{fooladivanda2017energy}; 2) the powers consumed by  the water pumps are then used an inputs to the AC optimal power flow problem (i.e., they are uncontrollable loads) solved by the power distribution operator, where the power curtailment from the V systems is minimized  along with the discrepancy from a given setpoint for the power at the substation~\cite{DallAnese17}. Table \ref{table:results} summarizes the achieved costs in the two cases. The total cost across the two systems is significantly lower in the OWPF case. It is clear that a coupled optimization approach enhances flexibility, which reduces the total cost compared to the decoupled optimization strategy. The operating cost of the water system is increased, and an increased amount of water in the tank at the end of the $ 12 $-hour optimization slot is observed; this calls for a systematic payment to the water system operator for providing services -- an interesting future research topic. On the other hand, the cost of PV power curtailment is almost halved, thus enabling a significantly higher utilization of renewable-based power. Fig.~\ref{fig:OWPF} illustrates the power profiles using  OWPF, while Fig.~\ref{fig:decoupled} shows the same profiles for the decoupled approach.

		\begin{table}[htbp]
			\centering
			\caption{Results of OWPF formulation and the decoupled approach.}
			\begin{tabular}{|c|c|c|}
				\hline
				{\textbf{Cost}} & {\textbf{OWPF}} &{\textbf{Decoupled Solver}} \\
				\hline
				\textbf{Water Cost}     & 460.26     & 268.47 \\
				\hline
				\textbf{Curtailment Cost}     & 769.30    & 1506.17 \\
				\hline
				\textbf{Substation Cost}     & 170.11   & 23.94 \\
				\hline
				\textbf{Total Cost}     &  1399.67  & 1798.58 \\
				\hline
			\end{tabular}%
			\label{table:results}%
		\end{table}%
		
		
		\begin{figure}[t]
			\centering
			\includegraphics[width=9.5cm]{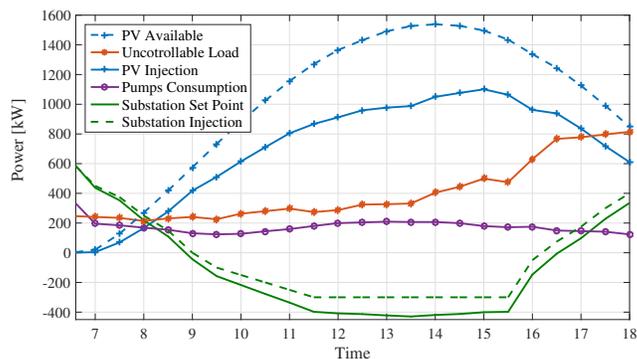}
			\caption{OWPF approach solution.}
			\label{fig:OWPF}
		\end{figure}
		
		\begin{figure}[t]
			\centering
			\includegraphics[width=9.5cm]{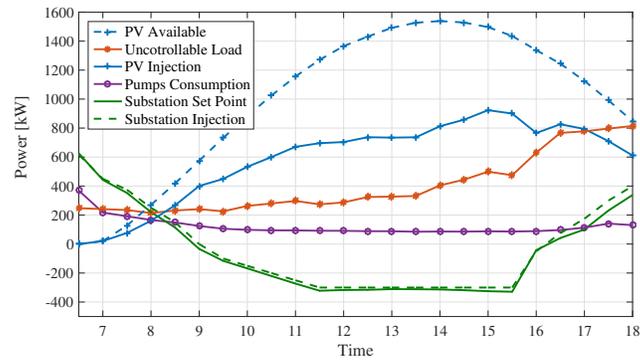}
			\caption{Decoupled approach solution.}
			\label{fig:decoupled}
		\end{figure}

		Finally, the convergence characteristics of the proposed distributed algorithms  are demonstrated using the same power and water networks. Each network operator initializes its local variables by optimizing a local operational cost without considering the other network constraints (i.e., the decoupled solution). Then, the distributed algorithm is used to reach a consensus on the powers consumed by the water pumps. The value of the parameter $ \rho $ is chosen to be $ 10^6 $. Fig.~\ref{fig:DOWPF-cost} shows the cost of operating the pumps at each local solver {per each iteration}. The discrepancy between the variables $ \beta_{ij}^t $ and $ p_{\sigma(ij)}^t $ is  shown in Fig.~\ref{fig:DOWPF-conv} by considering the Frobenius norm square measure of the difference between the two quantities in per unit. After approximately 50 iterations, the consensus error is lower than 1\%.

		\begin{figure}[t]
			\centering
			\includegraphics[width=9.5cm]{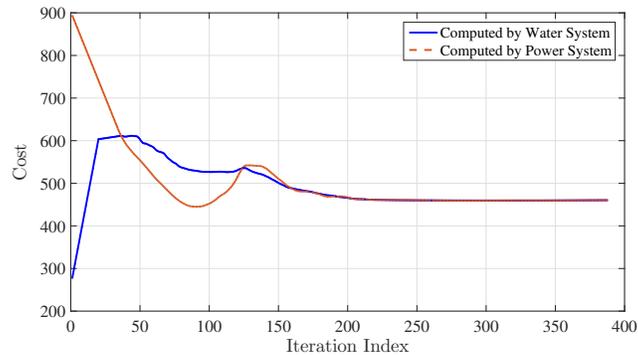}
			\caption{Operational cost of the water network by each system.}
			\label{fig:DOWPF-cost}
		\end{figure}
		
		\begin{figure}[t]
			\centering
			\includegraphics[width=9.5cm]{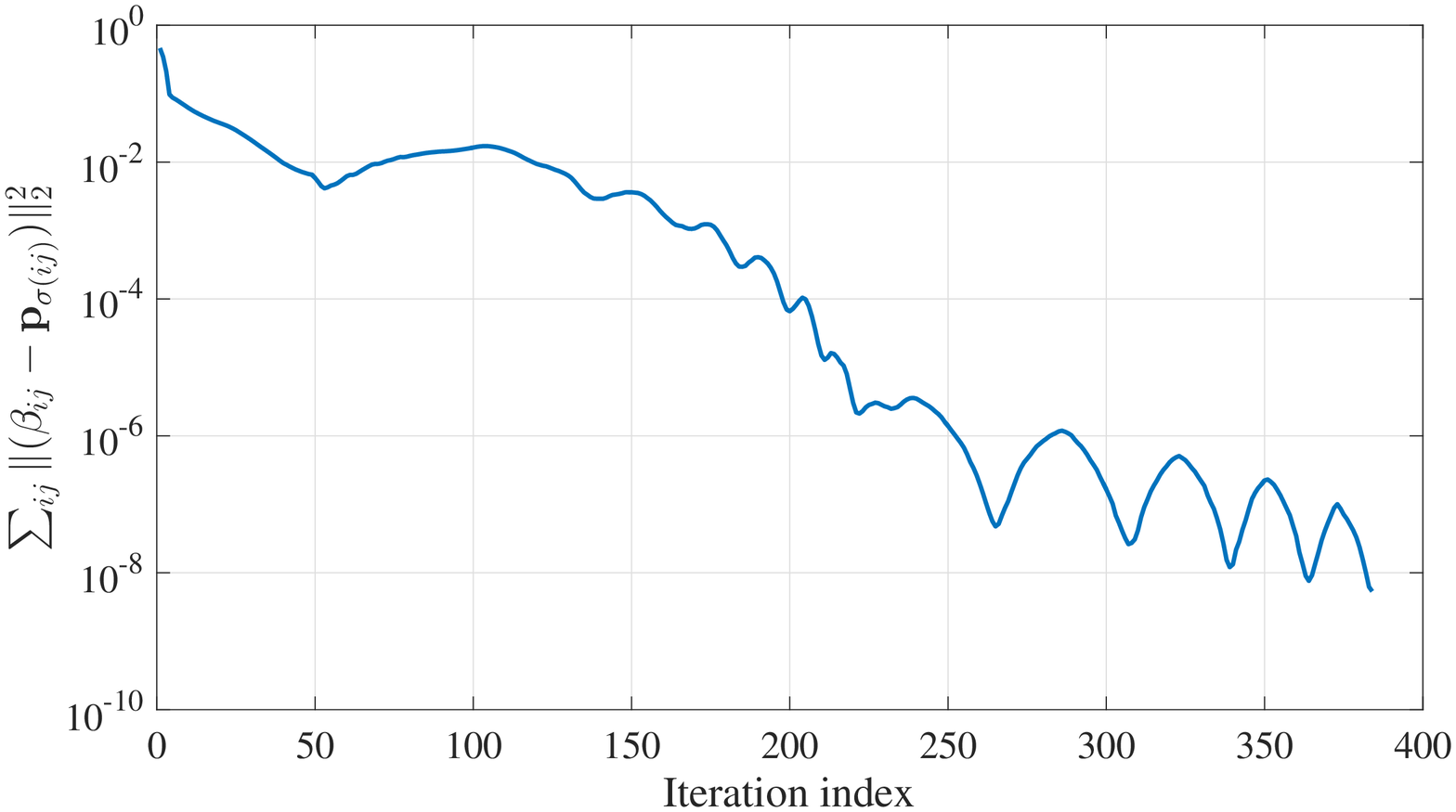}
			\caption{Convergence of the distributed solver.}
			\label{fig:DOWPF-conv}
		\end{figure}

		\section{Conclusions}
		This paper presented  an OWPF problem to optimize the use of controllable assets across power and water systems while accounting for the couplings between the two infrastructures. Although the physics governing the operation of the two systems and coupling constraints lead to a nonconvex  problem, feasible point pursuit-successive convex approximation approach was proposed  to identify feasible and optimal solutions. A distributed solver based on the ADMM was developed to enable water and power operators to pursue individual objectives while respecting the couplings between the two networks. The merits of the proposed approach were demonstrated for the case of a distribution feeder coupled with a municipal water distribution network. Future efforts will look at incorporating pump selection strategies into the OWPF, as well as the formulation of stochastic OWPF problems to account for  errors in the forecasts of power and water demands as well as powers available from renewable-based distributed energy resources.

		\bibliographystyle{IEEEtran}
		\bibliography{references}

\end{document}